\documentstyle[10pt]{article}
\begin{document}

\title{\bf \small THE RIEMANN HYPOTHESIS}
\author{\small Tribikram Pati}

\date{}

\maketitle
\begin{abstract} 
The Riemann Hypothesis is a conjecture made in 1859 by the great mathematician Riemann that
all the complex zeros of the zeta function $\zeta(s)$ lie on the `critical line' $\mbox{Rl}\,s=
1/2$. Our analysis shows that the assumption of the truth of the Riemann
Hypothesis leads to a contradiction. We are therefore led to the
conclusion that the Riemann Hypothesis is not true. 

\end{abstract}
\bigskip
\noindent 
1. The Zeta function of Riemann\footnote{For the theory of the Riemann Zeta 
Function, see Titchmarsh \cite{titchmarsh1}, Edwards \cite{edwards1}, Ivic \cite{ivic1}, 
\cite{ivic2}, Heath-Brown \cite{heathbrown1}, Bombieri \cite{bombieri1}, \cite{bombieri2}, 
 Montgomery \cite{montgomery1} and 
Conrey \cite{conrey1}, \cite{conrey2}. Also see Chandrasekharan \cite{chandrasekharan1}, 
Ramachandra \cite{ramachandra1} and Selberg \cite{selberg1}, \cite{selberg2}} 
is the analytic function obtained by the analytic continuation of the
sum-function $\zeta(s)$ of the infinite series

\begin{eqnarray}
\sum_{n=1}^{\infty} \frac{1}{n^{s}} \;\;\;\;(s=\sigma + it)
\label{eq1}
\end{eqnarray}
which is uniformly convergent in any finite region of the $s$-plane,
defined by $\sigma\geq 1+\delta, \delta > 0$. The Riemann Zeta function
$\zeta(s)$ is analytic in the entire $s$-plane, except for a simple pole
at $s=1$ with residue 1.

The function $\zeta(s)$ defined as the sum-function of (\ref{eq1}) is
also known from the researches of Euler to be defined as 
\begin{eqnarray}
\zeta(s)=\prod_{p}\left(1-\frac{1}{p^{s}}\right)^{-1},\;\;\;\;\sigma >
1,
\label{eq2}
\end{eqnarray}
the infinite product being taken for all prime numbers $p$. 

\bigskip

The Riemann Zeta function is given by
\begin{eqnarray} 
\zeta(s)= \frac{e^{-i\pi s}\Gamma(1-s)}{2\pi i}\int_{C} \frac{z^{s-
1}}{e^{z}-1}dz, 
\label{eq3}
\end{eqnarray}
where the contour $C$ consists of the real axis from $\infty$ to $\rho$
($0<\rho<2\pi$), the circle $\mid\!z\!\mid=\rho$ in the positive sense
round the origin and the real axis from $\rho$ to $\infty$, and $z^{s-
1}=\exp[(s-1)\log z]$, $\log z$ being assumed real at the starting point
of $C$. 

\bigskip

It follows from (\ref{eq2}) that $\zeta(s)\neq 0$ when $\sigma > 1$. It
has been proved independently by Hadamard \cite{hadamard1} and de la
Vall$\acute{\mbox{e}}$e-Poussin \cite{delavallee1} that $\zeta(s)\neq 0$ on the
line defined by $\sigma=1$. It is known that the only zeros of
$\zeta(s)$ on the $\sigma$ axis are at $s=-2, -4, -6,...$, and that all
the complex zeros of $\zeta(s)$ lie within the `critical strip' defined
by $0<\sigma < 1$ (See Titchmarsh \cite{titchmarsh1}). 

\bigskip

In his famous memoir ``$\ddot{\mbox{U}}$ber die Anzahl der Primzahlen unter einer
gegebenen Grosse'' of 1859 (for English translation, see Edwards
\cite{edwards1}) the celebrated German mathematician Riemann has
remarked that it is ``very likely'' that all the complex zeros of
$\zeta(s)$ have real part equal to $1/2$. This is called the `Riemann
Hypothesis'. This is equivalent to the statement that all the complex
zeros of $\zeta(s)$ lie on the `critical line' $\mbox{Rl}\,s
=\frac{1}{2}$. 

\bigskip

In 1914 the distinguished British mathematician G. H. Hardy
\cite{hardy1} proved that an infinity of complex zeros of $\zeta(s)$ lie
on the critical line. Subsequently Hardy and Littlewood \cite{hardy1},
Levinson \cite{levinson1}, Selberg \cite{selberg1} and Conrey
\cite{conrey1} have  estimated the proportion of complex zeros of
$\zeta(s)$ on the critical line to the number of complex zeros in the
critical strip. Bohr and Landau \cite{bohrlandau1} have also shown that
in the arbitrarily thin strip defined by $\frac{1}{2}-\epsilon < \sigma
<\frac{1}{2}+\epsilon$, $\epsilon > 0$, `almost all' the complex zeros of
$\zeta(s)$ lie. With the help of high-speed computers more and more
zeros of $\zeta(s)$ have been found to conform to Riemann's conjecture.
However, the Riemann Hypothesis has so far been neither proved nor
disproved. 

\bigskip

For an account of recent developments concerning the Riemann Hypothesis,
we may also refer to Marcus du Sautoy \cite{dusautoy1}, Dan Rockmore
\cite{rockmore1} and John Friedlander \cite{friedlander1}. 

\bigskip
\noindent
2. Let $\mu(n)$ be the M$\ddot{\mbox{o}}$bius function defined as:
$$
\mu(1)=1,\;\;\;\mu(n)=(-1)^{k},
$$
where $n$ is the product of $k$ different primes and $\mu(n)=0$ if $n$
contains as factor any prime raised to power $2$ of greater than $2$.
Let us write
$$
M(x)=\sum_{n\leq x}\mu(n).
$$
It was conjectured in 1897 by F. Mertens from numerical evidence that 
$$
M(n) < \sqrt n\;\;\;\;\;\;\;(n>1).
$$
This conjecture, called the `Mertens Hypothesis' was known to imply the
Riemann Hypothesis. It was proved to be false in October, 1983 by Andrew
Odlyzko and Hermann J. J. te Riele \cite{odlyzko1} after eight years of
historic collaboration, with the help of both mathematical techniques
and high-powered computation. 

\bigskip

3. In 1912, Littlewood \cite{littlewood1} proved 

\bigskip
\noindent
{\bf Theorem A.} The Riemann Hypothesis is equivalent to the statement
that, for $\epsilon > 0$, $M(x)= o\left(x^{\frac{1}{2}+\epsilon}\right)$,
as $x\rightarrow \infty$. 

\bigskip

In 1948, Turan \cite{turan1} proved that a sufficient condition for the
truth of the Riemann Hypothesis is that the functions 

\begin{eqnarray}
S_{n}(s)= \sum_{\nu =1}^{n} \frac{1}{\nu^{s}} \;\;\;\;(s=\sigma + it)
\label{eq4}
\end{eqnarray}
$n=1,2,3,...,$ should have no zeros in the half-plane $\sigma > 1$. He
had shown that this condition (\ref{eq4}) implies that, for all $x >0$, 
\begin{eqnarray}
\sum_{n\leq x} \frac{\lambda(n)}{n} \geq 0
\label{eq5}
\end{eqnarray}
where $\lambda(n)$ is the Liouville function, defined by: $\lambda(n)=(-
1)^{r}$, where $n$ has $r$ prime factors, a factor of multiplicity $k$
being counted $k$ times, and that (\ref{eq5}) implies the truth of the
Riemann Hypothesis. However, it was proved in 1958 by Haselgrove
\cite{haselgrove1} that (\ref{eq5}) is false in general. In 1968, Spira 
\cite{spira1} found that $S_{19}(s)$ has a zero in the region $\sigma >
1$. 

\bigskip

In 1885, T. J. Stieltjes \cite{stieltjes1} wrote to C. Hermite that he had
succeeded in proving that, as $x\rightarrow\infty$, $M(x)=
O(x^{\frac{1}{2}})$, and that this implies the truth of the Riemann
Hypothesis. As observed by Edwards \cite{edwards1}, in later years Stieltjes was 
unable to  reconstruct his proof. All he says about it is that it was
very difficult, that it was based on arithmetic arguments concerning
$\mu(n)$, and that he put it aside ``hoping to find a simpler proof of
the Riemann Hypothesis, based on the theory of the Zeta function rather
than on arithmetic.'' 

\bigskip
\noindent
4. The object of this paper is to prove the following.

\bigskip
\noindent
{\bf Theorem 1}. {\em The Riemann Hypothesis is not true. In other
words, not all the complex zeros of $\zeta(s)$ lie on the `critical
line':  $\sigma=1/2$ ($\sigma= \mbox{Rl}\,s$).}

\bigskip 
\noindent
{\bf Theorem 2}. {\em The `weakened Mertens Hypothesis': For every
$\epsilon > 0$, 
$$
M(x)=o(x^{\frac{1}{2}+\epsilon}),\;\;\;\mbox{as}\;\;x\rightarrow \infty,
$$
is not true.}

\bigskip
\noindent
{\bf Remarks}. In view of Theorem A, Theorem 2 is an immediate
consequence of Theorem 1, which also shows that Turan's condition stated
as (\ref{eq4}) in \S 3 is not satisfied. Some other immediate
consequences of Theorem 1 have been stated in the last section entitled
`Further Remarks'. 

\bigskip
\noindent
5.1 {\bf Lemmas.} For the proof of Theorem 1 we shall need a  number of
lemmas. We give these below. 

\bigskip
\noindent
{\bf N.B.} We shall sometimes write `RH' for `Riemann Hypothesis'. It is
also proper to state here that the assumption of RH is not necssarily
essential, but makes the proof of certain lemmas easier, and since we
assume the truth of the RH at the very outset of the proof of Theorem 1,
we assume it, wherever convenient, in the lemmas. However, some lemmas
are independent of the RH. 

\bigskip
\noindent
{\bf Lemma 1.} {\em Assume the RH. Then each interval $\left(\frac{1}{2}+in,
\frac{1}{2}+i(n+1)\right)$, $n$ large, of the `critical line' contains
a point $\frac{1}{2}+iT$, say, such that $\mid T - \gamma\mid >
\frac{A}{\log\,n}$, where $\gamma$ is the ordinate of any zero of
$\zeta(s)$ and $A$ is a positive absolute constant. Thus, the interval
of the `critical line' : $\left(\frac{1}{2}+i(T-A/\log n),
\frac{1}{2}+i(T + A/\log n)\right)$ is free from zeros of $\zeta(s)$} (Titchmarsh \cite{titchmarsh1}, p.340.)

\bigskip
\noindent
{\bf Proof:} It is known that the number $N$ of zeros of $\zeta(s)$ in the
interval $\left(\frac{1}{2}+in,\frac{1}{2}+i(n+1)\right)$ of the
critical line is given by the formula
$$
N \leq A_{1}\log\,n, \;\;\;\;\mbox{for}\;\;n\geq n_{0},
$$
where $n_{0}$ is sufficiently large, assuming as we may that $A_{1}$ is an absolute constant 
$\geq 1$.  Now if the $N$ zeros of $\zeta(s)$
are denoted by $z_{1}, z_{2},..,z_{N}$, then, assuming, if possible,
that the Lemma 1 is false, there is no point $T$ of the type described.
Then certainly the midpoints of the segments: $(\frac{1}{2}+in, z_{1}),
(z_{1},z_{2}),...,(z_{N}, \frac{1}{2}+i(n+1))$ are such that their
distances from the endpoints of the segments are each $\leq A/\log\,n$,
whatever positive constant $A$ may be, so that the total length of the segment, viz. 1,
satisfies 
\begin{eqnarray*}
1 &\leq & \frac{2A}{\log\,n}\times \mbox{no. of segments}\\
&= & \frac{2A}{\log\,n}\times (N+1)\leq \frac{2A}{\log\,n}(A_{1}\log
n +1)\\
& < &\frac{2A}{\log\,n}. 2A_{1}\log n\\
&=& 4AA_{1}.
\end{eqnarray*}

Hence $4AA_{1} > 1$ or $A> \frac{1}{4A_{1}}$. Taking
$A=\frac{1}{4A_{1}}$, we get a contradiction. Hence the Lemma is true. 

\bigskip
\noindent
{\bf Lemma 2.} {\em Assume the RH. Let $s=\sigma+it;\;
\;s_{0}=\frac{1}{2}+it_{0}$ is a zero of $\zeta(s)$ of order $m$ such
that $\mid t - t_{0}\mid < 1/\log\log\,t$, where $t, t_{0} > \tau$, a
sufficiently large positive number. Then, uniformly for $1/2\leq
\sigma\leq 2$, 
$$
\mbox{Rl}[\log \zeta(s) - m\log(s-s_{0})] < A\log t\frac{\log\log\log
t}{\log\log t},
$$
where $A$ is an absolute positive constant.}

\bigskip
\noindent
{\bf Proof.} The following theorem is known. 

\bigskip
\noindent
{\bf Theorem B.}\footnote{Titchmarsh \cite{titchmarsh1}, Th. 14.15}
Under the Riemann Hypothesis, uniformly for $1/2\leq\sigma\leq 2$, with
$s=\sigma + it$, as $t\rightarrow \infty$, 
$$
\log\zeta(s)-\sum_{\mid t -\gamma\mid <1/\log_{2}t}\log(s-\rho)=
O\left(\frac{\log t\log_{3}t}{\log_{2}t}\right),
$$
where $\rho=1/2+i\gamma$ runs through the zeros of $\zeta(s)$ and $\log_{2}t=
\log\log t$, $\log_{3}t=\log\log_{2}t$.  

\bigskip

By Theorem B, uniformly for $1/2\leq \sigma\leq 2$, 
$$
\log\frac{\zeta(s)}{\prod_{\mid t -\gamma\mid < 1/\log_{2}t}(s-\rho)} =
O\left(\widehat{\Omega}(t)\log_{3}t\right), 
$$
where $\widehat{\Omega}(t)\equiv \log t/\log_{2}t$. 

\bigskip

Hence, for $t>\tau$, where $\tau$ is a sufficiently large positive
number, 
\begin{eqnarray*}
\left| \log\left|\frac{\zeta(s)}{\prod_{\mid t -\gamma\mid <
1/\log_{2}t}(s-\rho)}\right|\right| &=& \left| 
\mbox{Rl}\,\log\frac{\zeta(s)}{\prod_{\mid t -\gamma\mid <1/\log_{2}t}(s-\rho)}
\right|\\
&\leq & \left| \log\frac{\zeta(s)}{\prod_{\mid t -\gamma\mid <
1/\log_{2}t}(s-\rho)}\right|< A\widehat{\Omega}(t)\log_{3}t. 
\end{eqnarray*}

Hence 
\begin{eqnarray}
\left| \frac{\zeta(s)}{\prod_{\mid t -\gamma\mid <
1/\log_{2}t}(s-\rho)}\right| < \exp (A\widehat{\Omega}(t)\log_{3}t), 
\label{eq6}
\end{eqnarray}

and 

\begin{eqnarray}
\left| \frac{\zeta(s)}{(s-s_{0})^{m}}\right| = 
\left| \frac{\zeta(s)}{\prod_{\mid t -\gamma\mid <
1/\log_{2}t}(s-\rho)}\right|\prod_{\mid t -\gamma\mid <
1/\log_{2}t,\rho\neq s_{0}} \mid s - \rho\mid.
\label{eq7}
\end{eqnarray}

The number of factors in $\prod_{\mid t -\gamma\mid <1/\log_{2}t}\mid s
-\rho\mid$, taking into account the multiplicity of the zeros, does not
exceed 
$$
N(t+\widehat{\delta}) - N(t-
\widehat{\delta})\;\;\;\;(\widehat{\delta}=1/\log_{2}t)
$$
which by Theorem 14.13 of Titchmarsh \cite{titchmarsh1}, 
$$
= O(\widehat{\delta}\log t)+ O(\widehat{\Omega}(t))=
O(\widehat{\Omega}(t)).
$$

Thus the number of factors in the product $\prod_{|t-
\gamma|<\widehat{\delta},\rho\neq s_{0}}\mid s -\rho\mid $ {\em a
fortiori} does not exceed $\widehat{A}\widehat{\Omega}(t)$. Also, since 
$$
s -\rho = \sigma - \frac{1}{2}+ i(t-\gamma)
$$
and $\mid t -\gamma\mid < \widehat{\delta}$, we have $\mid s -\rho\mid <
A$ (assuming, as we may, that $A$ is a positive absolute constant such that $\log A > 0$). Hence from
(\ref{eq6}) and (\ref{eq7}) we get 
\begin{eqnarray*}
\left| \frac{\zeta(s)}{(s-s_{0})^{m}}\right| &<&
\exp(A\widehat{\Omega}(t)\log_{3}t)A^{\widehat{A}\widehat{\Omega}(t)}\\
&=& \exp(A\widehat{\Omega}(t)\log_{3}t)\exp((\log
A)\widehat{A}\widehat{\Omega}(t))\\
&<& \exp(A\widehat{\Omega}(t)\log_{3}t).
\end{eqnarray*}

Hence, uniformly for $1/2\leq \sigma\leq 2$ and for $t,\;t_{0}>\tau$, where
$\tau$ is sufficiently large and positive, 
$$
\log \left| \frac{\zeta(s)}{(s-s_{0})^{m}}\right| <
A\widehat{\Omega}(t)\log_{3}t, 
$$
i.e. 
$$
\mbox{Rl}[\log\zeta(s) -\log (s-s_{0})^{m}] < A \frac{\log
t}{\log_{2}t}\log_{3}t, 
$$
whence the lemma follows. 

\bigskip
\noindent
{\bf Lemma 3.} {\em Assume the RH. If $s'=\sigma' + i\delta$ ($\sigma' >
1, \;\delta > 0$) is sufficiently close to the point $s'=1$, i.e. $ 0 <
\sigma' - 1 < \delta_{1},\;0<\delta <\delta_{1}$, where $\delta_{1}$ is
sufficiently small, then 
$$
\mbox{Rl}\log(s'-1) + \mbox{Rl}\log\zeta(s') < A,
$$
where $A$ is an absolute positive constant. }

\bigskip
\noindent
{\bf Proof:} When $s'$ is close to the point $s'=1$, then\footnote{
Titchmarsh \cite{titchmarsh1}, (2.1.16); Edwards \cite{edwards1}, p. 70} 
$$
\zeta(s')=\frac{1}{s'-1} + \gamma + O(\mid s' -1\mid).
$$
Hence 
$$
(s'-1)\zeta(s')= 1 + \gamma(s' -1) + O(\mid s' - 1\mid^{2}).
$$
Thus
$$
\lim_{s'\rightarrow 1}\left[(s'-1)\zeta(s')\right]=1.
$$
Hence 
$$
\lim_{s'\rightarrow 1}\left| (s'-1)\zeta(s')\right| = 1
$$
and, therefore, 
$$
\mid (s' -1)\zeta(s')\mid < 1+\epsilon
$$
for any $\epsilon > 0$, however small, if $0<\mid \sigma' -1\mid
<\delta_{1}$, $0<\delta <\delta_{1}$, where $\delta_{1}$ is sufficiently
small. Hence 
$$
\log\mid s' -1\mid +\log \mid \zeta(s')\mid <
\log(1+\epsilon)\;\;\;(\epsilon > 0),
$$
which gives the result. 

\bigskip
\noindent
{\bf Lemma 4.}\footnote{Titchmarsh \cite{titchmarsh1}, p. 46, (3.2.1)
and p. 189; Edwards \cite{edwards1}, p. 70} {\em For
$s'=\sigma'+it,\;\;\sigma' > 1$, 
$$
\log\zeta(s')=\sum_{p} \frac{1}{p^{s'}} + f(s')
$$
where 
$$
f(s')\equiv\sum_{p}\sum_{m=2}^{\infty} \frac{1}{mp^{ms'}}
$$
}
{\em is regular for $\sigma' >\frac{1}{2}$. Hence, for $\sigma' \geq 1$, 
$$
\mid f(s')\mid \leq \sum_{p}\sum_{m=2}^{\infty}
\frac{1}{p^{m}}=\sum_{p}\frac{1}{p(p-1)}=B, 
$$
say, a positive absolute
constant. Thus, since $\mid \mbox{Rl}f(s')\mid \leq \mid f(s')\mid$, 
$$
-B\leq \mbox{Rl} f(s')\leq B.
$$
} 

\bigskip
\noindent
{\bf Lemma 5.}\footnote{Titchmarsh \cite{titchmarsh1}, p. 46,
(3.2.2)}{\em If $\sigma'>1$, then, as $\sigma' \rightarrow 1$,} 
$$
\sum_{p} \frac{1}{p^{\sigma'}} \sim \log \frac{1}{\sigma' -1}.
$$

\bigskip
\noindent
{\bf Lemma 6.} {\em For $\delta > 0$ and any prime number $p$, 
$\left|\sin\left(\frac{\delta}{2}\log\,p\right)\right| \leq
\frac{\delta}{2}\log\,p.$}

\bigskip
\noindent
{\bf Proof.} If $0< \frac{\delta}{2}\log\,p\leq \frac{\pi}{2}$, 
$\left|\sin\left(\frac{\delta}{2}\log\,p\right)\right| = 
\sin\left(\frac{\delta}{2}\log\,p\right)\leq\frac{\delta}{2}\log\,p$ (see
Copson \cite{copson1}, p. 136). If $\frac{\pi}{2}<\frac{\delta}{2}\log\,p$, 
$\left|\sin\left(\frac{\delta}{2}\log\,p\right)\right|\leq 1<\frac{\pi}{2}
<\frac{\delta}{2}\log\,p$. 

\bigskip
\noindent
{\bf Lemma 7.} {\em Assume the RH. Then for $\sigma'=\mbox{Rl}s' > 1$,} 
$$
\sum_{p}\frac{\log p}{p^{\sigma'}}=\frac{1}{\sigma' -1} + O(1) =
O(\frac{1}{\sigma' -1}),
$$
{\em if $\sigma' - 1$ is bounded}. 

\bigskip
\noindent
{\bf Proof.} It is known that $(s'-1)\zeta(s')$ is regular in all finite
regions of the $s'$-plane. Also, under the RH, $(s'-1)\zeta(s')$ has no
zero in the strip $1/2<\sigma' < 1$, while, without RH, $(s'-
1)\zeta(s')$ has no zero in $\sigma' \geq 1$, it being known that (see
Lemma 3) $\lim_{s'\rightarrow 1}\left[ (s'-1)\zeta(s')\right]=1$, and
hence $\mid (s'-1)\zeta (s')\mid > 1-\epsilon$ for arbitrarily small
$\epsilon > 0$ ($\epsilon < 1/2$), for $\mid s' -1\mid < \delta'$,
$\delta'$ sufficiently small, and also that $(s'-1)\zeta(s')\neq 0$ in
the half-plane $\sigma' > 1$ and, by virtue of the result of Hadamard
\cite{hadamard1} and de la Vall$\acute{\mbox{e}}$e-Poussin
\cite{delavallee1}, $\zeta(s')$ has no zero on $\sigma'=1$. Hence, in
the entire half-plane $\mbox{Rl}s' =\sigma' > 1/2$, the function 
$$
\log\left[(s'-1)\zeta(s')\right]
$$
is regular, and hence for $\sigma' > 1$, 
\begin{eqnarray}
\frac{1}{\sigma' -1} + \frac{\zeta'(\sigma')}{\zeta(\sigma')}= O(1)
\label{eq8}
\end{eqnarray}

Also, for $\sigma' > 1$, 
$$
\log \zeta(s')=\sum_{p}\frac{1}{p^{s'}} + f(s')
$$
where $f(s')$ is regular for $\sigma' > 1/2$ (see Lemma 4). Hence for
$\sigma' > 1$, 
\begin{eqnarray}
\frac{\zeta'(\sigma')}{\zeta(\sigma')} &=& \frac{d}{d\sigma'}\left(
\sum_{p} \frac{1}{p^{\sigma'}}\right) + f'(\sigma') \nonumber\\
&=& -\sum_{p}\frac{\log p}{p^{\sigma'}} + f'(\sigma') \nonumber\\
&=&  -\sum_{p}\frac{\log p}{p^{\sigma'}} + O(1).
\label{eq9}
\end{eqnarray}

From (\ref{eq8}) and (\ref{eq9}), we have 
\begin{eqnarray*}
\sum_{p}\frac{\log p}{p^{\sigma'}}&=& -
\frac{\zeta'(\sigma')}{\zeta(\sigma')} + O(1)\\
&=& \frac{1}{\sigma' - 1} - O(1) + O(1) = \frac{1}{\sigma' -1} + O(1)\\
&=& O\left(\frac{1}{\sigma' -1}\right),
\end{eqnarray*}
if $\sigma' -1$ is bounded. 

\bigskip
\noindent
{\bf Lemma 8.} {\em Assume the RH. Then, with $s=\sigma+ it$, $\zeta(\frac{1}{2}+it)\neq 0$, 
uniformly for 
$1/2 <\sigma \leq 2$,} 
$$
{\cal E}\equiv \log\frac{1}{\mid \zeta(\sigma + it)\mid} - 
\log\frac{1}{\mid \zeta(1/2 + it)\mid}= O\left(\log t\left(\frac{\sigma -
1/2}{\delta}\right)\right),
$$
{\em where $0 <\delta =\min\mid t -\gamma\mid$, where $\gamma$ runs through
the zeros of $\zeta(s)$ such that $\mid t-\gamma\mid < 1$.}

\bigskip
\noindent
{\bf Proof.} Uniformly in $-1\leq \sigma\leq 2$, as $t\rightarrow
\infty$, \footnote{Titchmarsh \cite{titchmarsh1}, p. 340; see also Titchmarsh 
\cite{titchmarsh1}, Theorem 9.6(A)} 
\begin{eqnarray*}
\frac{\zeta'(s)}{\zeta(s)} &=& O\left(\frac{\log t}{\delta}\right) +
O(\log t)\\
&=&O\left(\frac{\log t}{\delta}\right)
\end{eqnarray*}
where $\delta$ is as defined above, so that $0<\delta < 1$. 

\bigskip

Hence, uniformly for $1/2 <\sigma \leq 2$, 
\begin{eqnarray*}
\mid {\cal E}\mid &=& \mid\log\mid\zeta(\sigma+it)\mid -\log \mid\zeta(1/2
+ it)\mid\mid \\
&=& \mid \mbox{Rl}\log \zeta(\sigma+it) - \mbox{Rl}\log
\zeta(1/2+it)\mid \\
&\leq& \mid \log \zeta(\sigma+it) - \log
\zeta(1/2+it)\mid\\
&=& \left|\int_{1/2}^{\sigma}\frac{\zeta'(x+it)}{\zeta(x+it)}dx\right|\\
&\leq &
\int_{1/2}^{\sigma}\left|\frac{\zeta'(x+it)}{\zeta(x+it)}\right|dx\\
&=& O\left(\frac{\log t}{\delta}\int_{1/2}^{\sigma}dx\right)=
O\left(\frac{\log t}{\delta}\left(\sigma - \frac{1}{2}\right)\right). 
\end{eqnarray*}

\bigskip
\noindent
{\bf Lemma 9.}\footnote{Titchmarsh \cite{titchmarsh1}, Th. 14.14 (A)} {\em Assume the RH. Then, for all $t > \tau$, a
sufficiently large positive number, 
$$
\mid\zeta(1/2+it)\mid \leq \widehat{A} e^{A_{3}\frac{\log t}{\log_{2}t}},
$$
where $\widehat{A}$ and $A_{3}$ are absolute positive constants, such
that $\widehat{A}>1$.}

\bigskip
\noindent
{\bf Lemma 10.} {\em If $a(t)$ and $k(t)$ are positive for $t
>\tau_{0}$ (a positive constant) and $\theta$ is a positive constant,
$a(t)\rightarrow \infty$ and $k(t)\rightarrow \infty$, as $t\rightarrow
\infty$, such that} 
$$
\frac{a(t)}{k(t)}\rightarrow
\infty\;\;\;\;\;\mbox{as}\;\;t\rightarrow\infty,
$$
{\em and, for $t > \tau_{0}$, $b(t)= 1 + \frac{\theta}{k(t)},$
then, for all $t \geq \tau_{1}$ (a sufficiently large positive
constant $> \tau_{0}$)}
$$
a(t)+b(t) < a(t)b(t).
$$

\bigskip
\noindent
{\bf Proof.} It is enough to show that, for $a(t)> 1$,
$$
b(t) > \frac{a(t)}{a(t)-1}
$$
and hence enough to show that 
$$
1+ \frac{\theta}{k(t)} > \frac{1}{1-\frac{1}{a(t)}}
$$
or 
$$
1 - \frac{1}{a(t)} +\frac{\theta}{k(t)}\left(1-\frac{1}{a(t)}\right) > 1
$$
or 
$$
\frac{\theta}{k(t)}\left(1-\frac{1}{a(t)}\right) > \frac{1}{a(t)}
$$
or 
$$
\theta > \frac{k(t)}{a(t)-1}= \frac{1}{\frac{a(t)}{k(t)}-
\frac{1}{k(t)}}.
$$
Since $\frac{a(t)}{k(t)}\rightarrow \infty$ and $k(t)\rightarrow\infty$,
the result follows from the assumption that $\theta$ is a constant $>
0$. 

\bigskip
\noindent
5.2 {\bf Proof of Theorem 1.}
We assume that the Riemann Hypothesis is true. We use the notations:
$$
\log_{2}t=\log\log t;\;\;\log_{3}t=\log(\log_{2}t), \;\;\Omega(t)=
\frac{\log t\log_{3}t}{\log_{2}t}
$$
Throughout $A$ will denote an absolute positive constant. For the sake of clarity we may 
sometimes use notations like $A^{\circ}, A_{1}, A_{2}, A'$ etc. for
specific absolute constants. It is understood that the notation $A$ does
not necessarily indicate the same positive absolute constant. 

\bigskip

We consider the interval $(n, n+1)$, where the positive integer $n$ is
sufficiently large as per the requirements of our analysis. By Lemma 1,
there will exist a point $\frac{1}{2}+iT$ on the critical line such that
$T$ lies in the interval $(n, n+1)$ and the region of the critical strip
defined by the horizontal lines with ordinates: $\left(T-\frac{A^{\circ}}{\log
n}\right)$, $\left(T+ \frac{A^{\circ}}{\log n}\right)$ is free from zeros of the
Zeta function $\zeta(s)$, $A^{\circ}$ being an absolute positive constant.  Let 
$s_{0}= \frac{1}{2}+it_{0}$ be the first zero of $\zeta(s)$ to be
encountered just below $\frac{1}{2}+ i\left(T-\frac{A^{\circ}}{\log
n}\right)$. Let us write
\begin{eqnarray}
t= t_{0}+\delta
\label{eq10}
\end{eqnarray}
where $\delta$ is defined as 
\begin{eqnarray}
\delta=\frac{A^{\ast}}{\log(n+1)}
\label{eq11}
\end{eqnarray}
where $A^{\ast}$ is a suitably small positive number to be further
specified in the sequel such that
\begin{eqnarray}
A^{\ast} <A^{\circ}.
\label{eq12}
\end{eqnarray}
By (\ref{eq11}) and (\ref{eq12}),
\begin{eqnarray}
\delta < \frac{A^{\circ}}{\log n}.
\label{eq13}
\end{eqnarray}
Now by (\ref{eq10}), (\ref{eq11}), (\ref{eq12}) and (\ref{eq13}) and the
definition of $T$, we have
\begin{eqnarray}
t=t_{0}+\delta &<& t_{0} + \frac{A^{\circ}}{\log n} < T < n+1;\nonumber\\
\delta &=& \frac{A^{\ast}}{\log(n+1)} < \frac{A^{\circ}}{\log t} <
\frac{1}{\log_{2}t}
\label{eq14}
\end{eqnarray}
for sufficiently large $t$, taking $n$ sufficiently large. 

\bigskip

Let $\frac{1}{2}+it_{1}$ be the first zero of $\zeta(s)$ to be encountered
just above $\frac{1}{2}+ i\left(T+\frac{A^{\circ}}{\log n}\right)$. Then 
$$
\mid t -t_{1}\mid = t_{1}-t = t_{1} -T +T-t,
$$
so that 
$$
\mid t -t_{1}\mid > t_{1}- T > \frac{A^{\circ}}{\log n},
$$
while, by (\ref{eq13}),
$$
\mid t - t_{0}\mid = t - t_{0} =\delta < \frac{A^{\circ}}{\log n}\;.
$$

\bigskip

Hence 
\begin{eqnarray}
\delta =\min\mid t -\gamma\mid,
\label{eq15}
\end{eqnarray}
where $\gamma$ runs through the zeros of $\zeta(s)$ such that $\mid t-
\gamma \mid < 1$. 

\bigskip

By virtue of (\ref{eq14}) and Lemma 2, uniformly for $1/2\leq \sigma\leq
2$, writing $\Omega(t)=\frac{\log t\log_{3}t}{\log_{2}t}$, we have 
\begin{eqnarray}
\mbox{Rl}\log\zeta(\sigma + it) - m\mbox{Rl}\log \left(\sigma -
\frac{1}{2} + i(t-t_{0})\right) < A\Omega(t). 
\label{eq16}
\end{eqnarray}

\bigskip

Hence uniformly in $1/2\leq \sigma \leq 2$
$$
\log \mid \zeta(\sigma + it)\mid - m\mbox{Rl}\log\left(\sigma -
\frac{1}{2} + i\delta\right) < A\Omega(t)
$$
or 
$$
\log \mid \zeta(\sigma + it)\mid - m\mbox{Rl}\log\left(\sigma' -1 + 
i\delta\right) < A\Omega(t)
$$
where $\sigma' -1 = \sigma - \frac{1}{2}$. 

\bigskip

Using Lemma 3, since $\sigma-\frac{1}{2}=\sigma' -1$, we have, uniformly
in $\sigma: \;0< \sigma -\frac{1}{2} < \delta_{1}$, with $\delta
<\delta_{1}$, 

$$
m\mbox{Rl}\log\zeta(\sigma' + i\delta) < A\Omega(t) + mA +
\log\frac{1}{\mid\zeta(\sigma+it)\mid}
$$

or 

\begin{eqnarray*}
\mbox{Rl}\log\zeta(\sigma' + i\delta) &<& \frac{A}{m}\Omega(t) + A +
\frac{1}{m}\log\frac{1}{\mid\zeta(\sigma+it)\mid}\\
&<& A\Omega(t) + \frac{1}{m}\log\frac{1}{\mid\zeta(\sigma+it)\mid}.
\end{eqnarray*}

Hence by Lemma 4, and using its notation, 
$$
\mbox{Rl}\sum_{p} \frac{1}{p^{\sigma'+i\delta}} +
\mbox{Rl}f(\sigma'+i\delta) < A\Omega(t) + \frac{1}{m}\log\frac{1}{\mid
\zeta(\sigma+it)\mid}\,.
$$

Appealing again to Lemma 4, with $B$ having the same meaning as in Lemma
4, 
$$
\mbox{Rl}\sum_{p} \frac{1}{p^{\sigma'+i\delta}} < A\Omega(t) + B + 
\frac{1}{m}\log\frac{1}{\mid\zeta(\sigma+it)\mid}
$$
or 
$$
\sum_{p}\frac{\cos(\delta\log p)}{p^{\sigma'}}< A\Omega(t) +
\frac{1}{m}\log\frac{1}{\mid\zeta(\sigma+it)\mid}
$$
or 
$$
\sum_{p}\frac{1}{p^{\sigma'}}< A\Omega(t)+\frac{1}{m}\log\frac{1}{\mid\zeta(\sigma+it)\mid}
+ 2\sum_{p}\frac{\sin^{2}(\frac{\delta}{2}\log p)}{p^{\sigma'}}.
$$

By Lemma 5, for $\eta > 0$, arbitrarily small ($<\frac{1}{2}$), $\exists \;
\sigma_{0}'(\eta)>1$ such that, uniformly for $\sigma' :
\;1<\sigma'<\sigma_{0}'(\eta)$, i.e. uniformly for $\sigma:\;\frac{1}{2} <
\sigma <\sigma_{0}(\eta)$, $\sigma_{0}(\eta)=\sigma_{0}'(\eta)-\frac{1}{2}$, 
$$
\sum_{p}\frac{1}{p^{\sigma'}}> (1-\eta)\log\frac{1}{\sigma' -1}.
$$
Hence, uniformly for $0<\sigma-\frac{1}{2}< \sigma_{0}^{\ast}-
\frac{1}{2} = \min\left[\sigma_{0}(\eta)-\frac{1}{2},
\delta_{1}\right]$, with $0< \delta <\delta_{1}$, 
$$
(1-\eta)\log\frac{1}{\sigma - \frac{1}{2}}<A\Omega(t)+ 
\frac{1}{m}\log\frac{1}{\mid\zeta(\sigma+it)\mid}+ 2\sum_{p}
\frac{\mid\sin(\frac{\delta}{2}\log p)\mid}{p^{\sigma'}}
$$
since $\sin^{2}(\frac{\delta}{2}\log p)=\mid\sin(\frac{\delta}{2}\log
p)\mid^{2}\leq \mid\sin(\frac{\delta}{2}\log p)\mid.$

Hence, by Lemma 6, uniformly for $0<\sigma-\frac{1}{2}<
\sigma_{0}^{\ast}-\frac{1}{2}$, with $0<\delta < \delta_{1}$, 
$$
\log\frac{1}{\sigma - \frac{1}{2}}<\frac{A}{1-\eta}\Omega(t)+ 
\frac{1}{(1-\eta)m}\log\frac{1}{\mid\zeta(\sigma+it)\mid}+ \frac{\delta}{1-\eta}
\sum_{p}\frac{\log p}{p^{\sigma'}}
$$
so that, by Lemma 7, uniformly for $0<\sigma -
\frac{1}{2}<\sigma_{0}^{\ast}-\frac{1}{2}$, with $0<\delta<\delta_{1}$, 
\begin{eqnarray*}
\log\frac{1}{\sigma - \frac{1}{2}} &<& 2A\Omega(t) +
\widetilde{A}\log\frac{1}{\mid\zeta(\sigma+it)\mid} +
A'\frac{\delta}{\sigma-\frac{1}{2}}\\
& &\left(\widetilde{A}=\frac{1}{(1-\eta)m}\leq \frac{1}{1-
\eta}< 2\;\;\mbox{and} \;A'\;\mbox{is an absolute positive constant}\right)\\
=A_{1}\Omega(t) &+& \widetilde{A}\log\frac{1}{\mid\zeta(\frac{1}{2}+it)\mid} +
\widetilde{A}\left[\log\frac{1}{\mid\zeta(\sigma+it)\mid}-
\log\frac{1}{\mid\zeta(\frac{1}{2}+it)\mid}\right] + A'\frac{\delta}{\sigma-\frac{1}{2}}\\
\leq  A_{1}\Omega(t) &+& \widetilde{A}\log\frac{1}{\mid\zeta(\frac{1}{2}+it)\mid} +
2\left|\log\frac{1}{\mid\zeta(\sigma+it)\mid}-
\log\frac{1}{\mid\zeta(\frac{1}{2}+it)\mid}\right| +
A'\frac{\delta}{\sigma-\frac{1}{2}}.
\end{eqnarray*}

Hence, by virtue of (\ref{eq15}) and Lemma 8, for sufficiently large
$n$, and hence $t$, uniformly for $\sigma:\;0<\sigma-
\frac{1}{2}<\sigma_{0}^{\ast}-\frac{1}{2}$, with $0<\delta <
\delta_{1}$, 
\begin{eqnarray}
\log\left(\frac{\mid\zeta(\frac{1}{2} +it)\mid^{\widetilde{A}}}{\sigma-\frac{1}{2}}\right)
<A_{1}\Omega(t) +A_{2}\log t\frac{\sigma -\frac{1}{2}}{\delta} +
A'\frac{\delta}{\sigma-\frac{1}{2}}.
\label{eq17}
\end{eqnarray}

We write 
\begin{eqnarray}
\widehat{\sigma}-\frac{1}{2} = \frac{\left|
\zeta\left(\frac{1}{2}+it\right)\right|^{\widetilde{A}}}{\exp(A_{1}\Omega + A_{2}\log t
+\theta)}\;,
\label{eq18}
\end{eqnarray}
where $\theta$ is a positive constant; $\Omega=\Omega(t)$. We observe
that, for all sufficiently large $t$, 
$$
\left(\frac{A_{1}\Omega +\theta}{A_{1}\Omega}\right)^{2} <
\frac{A_{2}}{A'}\log t.
$$
Hence, for sufficiently large $n$ and hence $t$, 
\begin{eqnarray}
\log(n+1)\frac{A_{1}\Omega +\theta}{A_{1}\Omega}\left(\widehat{\sigma}-
\frac{1}{2}\right)
<\log(n+1)\frac{A_{2}}{A'}\frac{A_{1}\Omega}{A_{1}\Omega+\theta}\log
t\left(\widehat{\sigma}-\frac{1}{2}\right).
\label{eq19}
\end{eqnarray}

Therefore, there exists a positive number $A^{\ast}$ such that 
\begin{eqnarray}
\log(n+1)\frac{A_{1}\Omega +\theta}{A_{1}\Omega}\left(\widehat{\sigma}-
\frac{1}{2}\right)
<A^{\ast}<\log(n+1)\frac{A_{2}}{A'}\frac{A_{1}\Omega}{A_{1}\Omega+\theta}\log
t\left(\widehat{\sigma}-\frac{1}{2}\right).
\label{eq20}
\end{eqnarray}

For $n > 3$, $(n+1)< 2(n-1)$, so that 
\begin{eqnarray}
R &\equiv & \log(n+1)\frac{A_{2}}{A'}\frac{A_{1}\Omega}{A_{1}\Omega+\theta}\log
t\left(\widehat{\sigma}-\frac{1}{2}\right)\nonumber\\
&< & \log[2(n-1)]\frac{A_{2}}{A'}\frac{A_{1}\Omega}{A_{1}\Omega+\theta}\log
t\left(\widehat{\sigma}-\frac{1}{2}\right)\nonumber\\
&<& \log(2t)\frac{A_{2}}{A'}\log t\left(\widehat{\sigma}-
\frac{1}{2}\right)\;\;\;\;\;\;\;\;\;\;[n-1 < t]\nonumber\\
&<& 2\log t\frac{A_{2}}{A'}\log t\left(\widehat{\sigma}-
\frac{1}{2}\right)\;\;\;\;\;\;\;\;\;\;[2< n-1< t]\nonumber \\
&=& \frac{2A_{2}}{A'}(\log t)^{2}\left(\widehat{\sigma}-
\frac{1}{2}\right), 
\label{eq21}
\end{eqnarray}
since, for $t>n-1 > 2$, $\log(2t)=\log 2 +\log t < 2\log t$. Thus from
(\ref{eq18}) and Lemma 9, 
\begin{eqnarray}
R &<& \frac{2A_{2}\widehat{A}^{\widetilde{A}}\exp\left(2\log_{2}t +A_{3}\widetilde{A}
\frac{\log t}{\log_{2}t}\right)}{A'\exp(A_{2}\log t+A_{1}\Omega(t)+\theta)}\nonumber\\
&=& \frac{2A_{2}\widehat{A}^{\widetilde{A}}}{A'\exp\left(A_{2}\log t
+\frac{\log t}{\log_{2}t}\left(A_{1}\log_{3}t-A_{3}\widetilde{A}\right)
- 2\log_{2}t +\theta\right)}.
\label{eq22}
\end{eqnarray}

Now, since $\widetilde{A}=\frac{1}{(1-\eta)m}\leq \frac{1}{1-\eta} < 2$,
$$
A_{1}\log_{3}t - A_{3}\widetilde{A} > A_{1}\log_{3}t - 2A_{3},
$$
which is positive for all sufficiently large $t$. 

\bigskip

Hence from (\ref{eq22}), for sufficently large $n$, and hence $t$, 
\begin{eqnarray}
R < \frac{2A_{2}\widehat{A}^{2}}{A'\exp\left(A_{2}\log t
+\frac{\log t}{\log_{2}t}\left(A_{1}\log_{3}t-2A_{3}\right)
- 2\log_{2}t +\theta\right)}
< A^{\circ}
\label{eq23}
\end{eqnarray}
however small $A^{\circ}$ may be, since $\widehat{A} > 1$, as assumed in
Lemma 9. 

From (\ref{eq20}) and (\ref{eq23}) it follows that 
\begin{eqnarray}
A^{\ast} < A^{\circ},
\label{eq24}
\end{eqnarray}
which had been stated in (\ref{eq12}), as an assumption. 

\bigskip

From (\ref{eq18}) and (\ref{eq20}), it follows that for sufficiently
large $n$ and hence $t$, 

\begin{eqnarray}
\delta &=& \frac{A^{\ast}}{\log(n+1)} < \frac{A_{2}}{A'}\log
t\left(\widehat{\sigma}-\frac{1}{2}\right)\nonumber \\
&<& \frac{A_{2}\widehat{A}^{\widetilde{A}}}{A'\exp\left(
A_{2}\log t+\frac{\log t}{\log_{2}t}(A_{1}\log_{3}t - 2A_{3})-\log_{2}t+\theta\right)}\;\;\;\;\;
\left[\mbox{since} \widetilde{A}\leq \frac{1}{1-\eta}< 2\right]\nonumber\\
&<& \frac{A_{2}\widehat{A}^{2}}{A'\exp\left(
A_{2}\log t+\frac{\log t}{\log_{2}t}(A_{1}\log_{3}t - 2A_{3})-\log_{2}t+\theta\right)}
\nonumber\\
&<& \sigma_{0}^{\ast} - \frac{1}{2}
\label{eq25}
\end {eqnarray}
however small $\sigma_{0}^{\ast}-\frac{1}{2}$ may be. 

\bigskip

Now, from (\ref{eq17}) it follows that, for sufficiently large $n$, and
hence $t$, uniformly in $\sigma: \;0< \sigma-\frac{1}{2} < \delta <
\sigma_{0}^{\ast}-\frac{1}{2}$,

\begin{eqnarray}
\frac{1}{A'}\log\left(\frac{\mid
\zeta(\frac{1}{2}+it)\mid^{\widetilde{A}}}{\sigma-1/2}\right) < 
\frac{1}{A'}\left(A_{1}\Omega+ A_{2}\log t\frac{\sigma-
\frac{1}{2}}{\delta}\right)+\frac{\delta}{\sigma-1/2}.
\label{eq26}
\end{eqnarray}

\bigskip

Let us now take
\begin{eqnarray}
\overline{\sigma}-
\frac{1}{2}=\frac{A_{1}\Omega}{A_{1}\Omega+\theta}\delta.
\label{eq27}
\end{eqnarray}

\bigskip

Hence from (\ref{eq25}) 

\begin{eqnarray}
\overline{\sigma}-\frac{1}{2}<\delta < \sigma^{\ast}_{0}-\frac{1}{2}.
\label{eq28}
\end{eqnarray}

\bigskip

We use the notations:
\begin{eqnarray}
a(t)\equiv \frac{1}{A'}\left(A_{1}\Omega+A_{2}\log
t\frac{\overline{\sigma}-\frac{1}{2}}{\delta}\right),
\label{eq29}
\end{eqnarray}

\begin{eqnarray}
b(t)\equiv \frac{\delta}{\overline{\sigma}-1/2}=\frac{A_{1}\Omega+\theta}{A_{1}\Omega}= 1
+\frac{\theta}{A_{1}\Omega}=1+\frac{\theta}{k(t)},
\label{eq30}
\end{eqnarray}

where

\begin{eqnarray}
k(t)\equiv A_{1}\Omega(t).
\label{eq31}
\end{eqnarray}

Evidently, by (\ref{eq27}) and (\ref{eq29}), $a(t)\rightarrow\infty$, as
$t\rightarrow\infty$; also

\begin{eqnarray}
k(t)\rightarrow \infty, \;\;\;\mbox{as}\;\;t\rightarrow\infty.
\label{eq32}
\end{eqnarray}

Also, as $t\rightarrow\infty$, 
\begin{eqnarray}
\frac{a(t)}{k(t)} &=& \frac{1}{A'}\left(1+ \frac{A_{2}\log
t}{A_{1}\Omega}\frac{A_{1}\Omega}{A_{1}\Omega+\theta}\right)\nonumber\\
&=& \frac{1}{A'}\left(1+ \frac{A_{2}}{A_{1}}\frac{\log_{2}t}{\log_{3}t}
\frac{A_{1}\Omega}{A_{1}\Omega+\theta}\right)\rightarrow\infty
\label{eq33}
\end{eqnarray}

\bigskip

Therefore, by Lemma 10, 
$$
a(t)+b(t) < a(t).b(t).
$$
Thus from (\ref{eq26}), we get, for $\overline{\sigma}-\frac{1}{2}$ as
per (\ref{eq27}), using (\ref{eq25}),
\begin{eqnarray}
\log\left(\frac{\mid\zeta
\left(\frac{1}{2}+it\right)\mid^{\widetilde{A}}}{\overline{\sigma}-
1/2}\right) &<& \left(A_{1}\Omega + A_{2}\log t\frac{\overline{\sigma}-
\frac{1}{2}}{\delta}\right)\frac{\delta}{\overline{\sigma}-
1/2}\nonumber\\
&=& A_{1}\Omega\frac{\delta}{\overline{\sigma}-1/2}+ A_{2}\log
t\nonumber\\
&=& A_{1}\Omega \frac{A_{1}\Omega+\theta}{A_{1}\Omega} + A_{2}\log
t\nonumber\\
&=&A_{1}\Omega+A_{2}\log t+\theta
\label{eq34}
\end{eqnarray}
or
\begin{eqnarray}
\log\left(\frac{\mid\zeta\left(\frac{1}{2}+it\right)\mid^{\widetilde{A}}(A_{1}\Omega+
\theta)}{A_{1}\Omega.\delta}\right)< A_{1}\Omega+A_{2}\log t+\theta,
\label{eq35}
\end{eqnarray}

whence we have

\begin{eqnarray}
\delta > \frac{A_{1}\Omega +\theta}{A_{1}\Omega}
\frac{\mid\zeta\left(\frac{1}{2}+it\right)\mid^{\widetilde{A}}}{\exp(A_{1}\Omega+A_{2}\log
t+\theta)}
= \frac{A_{1}\Omega+\theta}{A_{1}\Omega}\left(\widehat{\sigma}-
\frac{1}{2}\right)
\label{eq36}
\end{eqnarray}
(which is equivalent to the first inequality in (\ref{eq20}) by the
definition of $\delta$ (see (\ref{eq11})) and (\ref{eq18}), the definition of
$\widehat{\sigma}-\frac{1}{2}$. 

\bigskip

If we take 
\begin{eqnarray}
\delta' = \frac{A_{1}\Omega +\theta}{A_{1}\Omega}\left(\widehat{\sigma}-
\frac{1}{2}\right),
\label{eq37}
\end{eqnarray}
then, from (\ref{eq36}) we have 
\begin{eqnarray}
\delta' <\delta.
\label{eq38}
\end{eqnarray}

Thus
\begin{eqnarray}
\widehat{\sigma}-\frac{1}{2}=
\frac{A_{1}\Omega}{A_{1}\Omega+\theta}\delta' < \frac{A_{1}\Omega}{A_{1}\Omega+\theta}\delta
<\delta
\label{eq39}
\end{eqnarray}

Since $\widehat{\sigma}-\frac{1}{2}<\delta$, it follows from
(\ref{eq25}) that 
\begin{eqnarray}
\widehat{\sigma}-\frac{1}{2}<\sigma_{0}^{\ast}-\frac{1}{2}
\label{eq40}
\end{eqnarray}
so that by virtue of (\ref{eq17}) and (\ref{eq25}), 
\begin{eqnarray}
\log\left(\frac{\mid\zeta
\left(\frac{1}{2}+it\right)\mid^{\widetilde{A}}}{\widehat{\sigma}-
1/2}\right)< A_{1}\Omega+ A_{2}\log t\frac{\widehat{\sigma}-
\frac{1}{2}}{\delta} + A'\frac{\delta}{\widehat{\sigma}-\frac{1}{2}}.
\label{eq41}
\end{eqnarray}
We show below that 
\begin{eqnarray} 
A_{2}\log t\frac{\widehat{\sigma}-
\frac{1}{2}}{\delta} + A'\frac{\delta}{\widehat{\sigma}-\frac{1}{2}}< 
A_{2}\log t\frac{\widehat{\sigma}-
\frac{1}{2}}{\delta'} + A'\frac{\delta'}{\widehat{\sigma}-\frac{1}{2}}.
\label{eq42}
\end{eqnarray}
The inequality (\ref{eq42}) holds if 
$$
A_{2}\log t\left(\widehat{\sigma}-
\frac{1}{2}\right)\left(\frac{1}{\delta'}-\frac{1}{\delta}\right) +
A'\frac{\delta' -\delta}{\widehat{\sigma}-\frac{1}{2}} > 0
$$
or 
$$
A_{2}\log t\left(\widehat{\sigma}-
\frac{1}{2}\right)^{2}- A'\delta'\delta > 0
$$
or, using (\ref{eq37}), if 
$$
A_{2}\log t\frac{A_{1}\Omega}{A_{1}\Omega+\theta}\left(\widehat{\sigma}-
\frac{1}{2}\right)- A'\delta > 0\;\;\;\;\;\;\;\;[(\ref{eq10}): \delta=A^{\ast}/\log(n+1)].
$$
or, by (\ref{eq10}), 
$$
A_{2}\log t\frac{A_{1}\Omega}{A_{1}\Omega+\theta}\left(\widehat{\sigma}-
\frac{1}{2}\right) > A'\frac{A^{\ast}}{\log(n+1)}
$$
or 
$$
A^{\ast} < \log(n+1)\frac{A_{2}}{A'}\log t\frac{A_{1}\Omega}{A_{1}\Omega +\theta}\left(
\widehat{\sigma}-\frac{1}{2}\right),
$$
which holds by virtue of (\ref{eq20}). 

\bigskip

Thus (\ref{eq42}) is proved, and from (\ref{eq41}) we have 
\begin{eqnarray}
\log\left(\frac{\mid\zeta
\left(\frac{1}{2}+it\right)\mid^{\widetilde{A}}}{\widehat{\sigma}-
1/2}\right)< A_{1}\Omega + A_{2}\log t\frac{\widehat{\sigma}-
\frac{1}{2}}{\delta'} + A'\frac{\delta'}{\widehat{\sigma}-\frac{1}{2}}.
\label{eq43}
\end{eqnarray}

Hence

\begin{eqnarray} 
\frac{1}{A'}\log\left(\frac{\mid\zeta
\left(\frac{1}{2}+it\right)\mid^{\widetilde{A}}}{\widehat{\sigma}-
1/2}\right)< \frac{1}{A'}\left(A_{1}\Omega + A_{2}\log t\frac{\widehat{\sigma}-
\frac{1}{2}}{\delta'}\right) + \frac{\delta'}{\widehat{\sigma}-\frac{1}{2}}.
\label{eq44}
\end{eqnarray}

We note that, by (\ref{eq37}), 
$$
\widehat{\sigma}-\frac{1}{2} =
\frac{A_{1}\Omega}{A_{1}\Omega+\theta}\delta'. 
$$
We take 
\begin{eqnarray*}
\widehat{a}(t) &\equiv& \frac{1}{A'}\left( A_{1}\Omega+ A_{2}\log
t\frac{\widehat{\sigma}-\frac{1}{2}}{\delta'}\right),\\
\widehat{b}(t) &\equiv& \frac{\delta'}{\widehat{\sigma}-\frac{1}{2}},\\
k(t) &=& A_{1}\Omega(t).
\end{eqnarray*}
We find that 
$$
\widehat{b}(t)=1+\frac{\theta}{k(t)},
$$
and 
$$
\widehat{a}(t)\rightarrow\infty, \;\;k(t) \rightarrow \infty,\;\;\;\;\mbox{as}\;\;\;t\rightarrow\infty
$$
and, as in (\ref{eq33}),
$$
\frac{\widehat{a}(t)}{k(t)}\rightarrow \infty,
\;\;\;\mbox{as}\;\;t\rightarrow \infty.
$$
Hence, applying Lemma 10, we have 
$$
\widehat{a}(t)+ \widehat{b}(t) < \widehat{a}(t)\widehat{b}(t).
$$
Hence, from (\ref{eq44}) we get 
\begin{eqnarray}
\log\left(\frac{\mid\zeta
\left(\frac{1}{2}+it\right)\mid^{\widetilde{A}}}{\widehat{\sigma}-
1/2}\right) < A_{1}\Omega \frac{\delta'}{\widehat{\sigma}-\frac{1}{2}} +
A_{2}\log t = A_{1}\Omega(t) + A_{2}\log t+\theta.
\label{eq45}
\end{eqnarray}

Using (\ref{eq18}), we get 
$$
\log e^{A_{1}\Omega+ A_{2}\log t+\theta} < A_{1}\Omega + A_{2}\log t
+\theta,
$$
which is obviously false. 

\bigskip

This contradiction shows that the assumption that the Riemann Hypothesis
is true is not tenable. Thus the proof of Theorem 1 is complete. 

\bigskip

6. {\bf Further Remarks} 
\bigskip
\noindent
(I) It is known \footnote{Edwards \cite{edwards1}, \S 5.5.} that the
Riemann Hypothesis is equivalent to the statement that, for any
$\epsilon > 0$, and all sufficiently large $x>0$, the relative error in
the Prime Number Theorem:
$$
\pi(x)\sim \mbox{li}x, \;\mbox{where li}x=\lim_{\epsilon\rightarrow 0}\left[ \int_{0}^{1-
\epsilon} \frac{dt}{\log t} + \int_{1+\epsilon}^{x} \frac{dt}{\log t}
\right],
$$ 
is less that $x^{-\frac{1}{2}+\epsilon}$. 

In view of our Theorem 1, this statement is false. 

\bigskip
\noindent
(II) It is known that\footnote{Titchmarsh \cite{titchmarsh1}, p. 370, Theorem
14.25(B)} the following assertion is equivalent to the Riemann
Hypothesis. 
The series
$$
\sum_{n=1}^{\infty}\frac{\mu(n)}{n^{s}}
$$
is convergent for $\sigma = \mbox{Rl}s > \frac{1}{2}$. In view of
Theorem 1, this assertion is false. 

\bigskip
\noindent
(III) The following necessary and sufficient conditions for the truth of
the Riemann Hypothesis were given respectively by Hardy and Littlewood 
\footnote{Hardy and Littlewood \cite{hardylittlewood1}} and M. Riesz
\footnote{M. Riesz \cite{riesz1}}. By Theorem 1, these are false. 

\bigskip
\noindent
\begin{eqnarray*}
(\alpha)\;\;\;\;\;&\sum_{k=1}^{\infty}& \frac{(-x)^{k}}{k\,!\zeta(2k+1)} =
O(x^{-\frac{1}{4}}),\;\;\mbox{as}\;\;x\rightarrow\infty;\\
(\beta)\;\;\;\;\; &\sum_{k=1}^{\infty}& \frac{(-1)^{k+1}x^{k}}{(k-1)\,!\zeta(2k)} =
O(x^{\frac{1}{4}+\epsilon}),\;\;\mbox{as}\;\;x\rightarrow\infty,
\end{eqnarray*}
for every $\epsilon > 0$.

\vspace{1.5in}

\begin{verse}
10, Bank Road,\\
Allahabad 211002,\\
India\\
E-mail: pati\_tribikram@yahoo.co.in\\
Telephone: 91 532 2471052. 

\end{verse}
\end{document}